\documentclass[a4paper]{amsart}
\usepackage{amsmath,amssymb,amsfonts}
\usepackage{mathrsfs,latexsym,amsthm,enumerate}
\usepackage{amscd}
\usepackage[all]{xy}
\usepackage{tikz}
\usetikzlibrary{matrix,arrows}

\newcommand{\LCBS}{\mathsf{BS}}

\newcommand{\GBA}{\mathsf{BA}}

\newcommand{\D}{\mathcal{D}}

\newcommand{\LSBA}{\mathsf{Skew}}

\newcommand{\ESLCBS}{\mathsf{Etale}}

\newcommand{\Sp}{\mathbf{S}}
\newcommand{\Al}{\mathbf{A}}

\newcommand{\Hom}{\mathrm{Hom}}
\newcommand{\ocap}{\overline{\cap}}
\newcommand{\ucup}{\underline{\cup}}

\newtheorem{theorem}{Theorem}
\newtheorem{lemma}[theorem]{Lemma}
\newtheorem{corollary}[theorem]{Corollary}

\newtheorem{prop  osition}[theorem]{Proposition}
\newtheorem{remark}[theorem]{Remark}

\title[A dualizing object approach to non-commutative Stone duality]{A dualizing object approach to non-commutative Stone duality}
\author{Ganna Kudryavtseva}\thanks{Author's research is partially supported by ARRS grant P1-0288.}
\address{G. K.: University of Ljubljana,
Faculty of Computer and Information Science, \newline
Tr\v{z}a\v{s}ka cesta 25,
SI-1001, Ljubljana,
SLOVENIA.}
\email{ganna.kudryavtseva\symbol{64}fri.uni-lj.si}

\begin{document}
\maketitle
\begin{abstract} The aim of the present paper is to extend the dualizing object approach to Stone duality to the non-commutative setting of skew Boolean algebras. This continues the study of non-commutative generalizations of different forms of Stone duality   initiated in recent papers by Bauer and Cvetko-Vah, Lawson, Lawson and Lenz, Resende, and also the current author. In particular, we construct a series of dual adjunctions between the categories of Boolean spaces and skew Boolean algebras, unital versions of which are induced by dualizing objects $\{0,1,\dots, n+1\}$, $n\geq 0$. We describe Eilenberg-Moore categories of the monads of our adjunctions and construct  easily understood non-commutative reflections of skew Boolean algebras, where the latter can be faithfully embedded  (if $n\geq 1$) in a canonical way. As an application, we answer the question that arose in a recent paper by Leech and Spinks to describe the left adjoint to their `twisted product' functor~$\omega$.

\end{abstract}
\vspace{0.2cm}


\vspace{0.2cm}


\section{Introduction}

Stone duality \cite{S,D,GH} is probably the most famous result about Boolean algebras. It provides a subtle link between algebra and topology and has far-reaching consequences and numerous important generalizations. 

Recently, different variations of Stone duality have been generalized to non-commutative setting at the algebraic side to skew Boolean algebras by the author \cite{K}, to skew Boolean algebras with intersections independently by Bauer and Cvetko-Vah \cite{BCV} and the author \cite{K}, to Boolean inverse semigroups  by Lawson \cite{Law,Law1},  to pseudogroups   by Lawson and Lenz \cite{LL} and to certain classes of quantales by  Resende \cite{R,R1}. Very recently, Mark Lawson and the author \cite{KL4} have found a common generalization of results of \cite{K} and \cite{Law,Law1} to a duality where \'{e}tale actions of Boolean inverse semigroups arise at the algebraic side.  The latter actions were introduced in \cite{FS, St} and play an important role in Morita theory of inverse semigroups \cite{FLS, St} as well as in interplay of semigroup and topos theories \cite{F, FH, FS}. 

The study of skew Boolean algebras, that are objects of consideration of this paper, was initiated by Cornish \cite{C} and Leech \cite{L1,L2}. In \cite{BL} it was observed that skew Boolean algebras with intersections form a discriminator variety and  therefore tools of universal algebra can be applied to their study.
Stone duality for skew Boolean algebras \cite{K} and skew Boolean algebras with intersections  \cite{K, BCV} opens completely new perspectives of looking at these algebras. The duality of \cite{K} made the constructions of this paper possible, and also, together with results of \cite{Law,Law1}, led Mark Lawson and the author to the duality of \cite{KL4}. The latter duality, in its turn, witnesses 
that  skew Boolean algebras and Boolean inverse semigroups are closely related.  Putting them together leads to the construction of new objects that are connected with important concepts of  inverse semigroup theory and its applications.

An important feature of Stone duality is that (the unital version of) it is induced by a {\em dualizing object} (sometimes also called a {\em schizophrenic object}) $\{0,1\}$,  considered as a Boolean algebra or as a discrete topological space (see Subsection \ref{s21} for more details). To the contrast, its generalizations to skew Boolean setting \cite{BCV, K} {\em do~not} possess this property (see Subsection \ref{s24} for more details). The goal of the present paper is  to extend  the dualizing object approach to Stone duality to the non-commutative setting of skew Boolean algebras.  Our motivation is well-established significance of dualizing objects, both from the universal algebra and the category theory viewpoints, (see, e.g., \cite{CD,J, PT}), on the one hand, and recently revealed importance of skew Boolean algebras, on the other hand.

The non-commutative dualizing objects that arise in this paper are primitive left-handed skew Boolean algebras ${\bf n+2}=\{0,1,\dots, n+1\}$, $n\geq 0$ (see Subsection~\ref{s22} for the definitions).  Their role in our theory is similar to that of the two-element Boolean algebra ${\bf 2}=\{0,1\}$ in  Stone duality.  These objects induce a series of dual adjunctions $\Lambda_n \dashv \lambda_n$, $n\geq 0$, between the categories  of Boolean spaces and  left-handed skew Boolean algebras.

Given a left-handed skew Boolean algebra $S$, the Boolean algebra $S/{\mathcal D}$ {\em reflects} $S$ in the sense that there is a functor from left-handed skew Boolean algebras to Boolean algebras sending $S$ to $S/{\mathcal D}$ that is the left adjoint to the  inclusion functor in the reverse direction. Thus $S/{\mathcal D}$ is known as a {\em commutative reflection} of $S$. In this paper we construct a series $\lambda_n\Lambda_n(S)$, $n\geq 0$, of  reflections of $S$. If $n=0$ then $\lambda_0\Lambda_0(S)\simeq S/{\mathcal D}$. If $n\geq 1$ these reflections are non-commutative and
the units of the constructed adjunctions provide a canonical way to faithfully represent  $S$ in its very `Boolean algebra like' {\em non-commutative reflection} $\lambda_n\Lambda_n(S)$. However, in order to decrease the `degree of non-commutativity' of the enveloping algebra one has to sacrifise the size of the underlying Boolean algebra. 
Note that the possibility of an embedding of $S$ into $\lambda_n(X)$ for some $X$  can be easily deduced from \cite[Corollary 3.6]{LS} and Remark \ref{rem:ls1} of this paper. But no specific construction of an embedding of a skew Boolean algebra into another skew Boolean algebra with `low degree of non-commutativity' had been known before.

It is interesting that a variation of the functor $\lambda_1$ of  this paper has previously appeared in another disguise in the paper by Leech and Spinks \cite{LS} as the `twisted product' $\omega$-functor. Invoking the Freyd's adjoint functor theorem, it was argued in \cite{LS} that the functor $\omega$ has a left adjoint, $\Omega$, and the question to describe $\Omega$ arose. As a first step,  the action of $\Omega$ on finite objects was described in \cite{LS}. The approach of \cite{LS} did not allow to go further than that. Our approach provides a natural interpretation for both $\omega$ and $\Omega$ and leads to the full  description of  the functor $\Omega$ (see Remarks \ref{rem:ls1} and \ref{rem:ls2}). 

We also study the monads induced by the adjunctions $\Lambda_n \dashv \lambda_n$, $n\geq 0$, and  provide their nice combinatorial description.  We then prove that the Eilenberg-Moore categories of these monads are equivalent to the image of $\lambda_n$, so that this image is a reflective subcategory of the category of left-handed skew Boolean algebras. 

Before stating in Sections \ref{s4}, \ref{s6} and \ref{s7} the main constructions and results of this paper, we collect in Section \ref{s2} all necessary preliminaries. In particular, we explain what precisely we mean by a dualizing object approach to the classical Stone duality, then provide necessary background on skew Boolean algebras and \'{e}tale spaces and  review the non-commutative Stone duality from \cite{K} needed in  this paper. We also explain why this duality is not induced by a dualizing object.

\section*{Acknowledgements} I would like to thank the anonymous referee for suggesting that the topology on the space $\Lambda_n(S)$ might coincide with the topology inherited by $\Lambda_n(S)$ from the product space $\{0,\dots, n+1\}^S$. This turned out to be the case and led to a significant simplification of  the proof that $\Lambda_n(S)$ is a Boolean space in Section \ref{s6}. I am also grateful to the referee, Marcel Jackson and Mark Lawson for helpful comments.

\section{Preliminaries}\label{s2}

\subsection{The dualizing object view of the classical Stone duality}\label{s21}
By a {\em Boolean algebra} we mean what is usually called a generalized Boolean algebra, that is a relatively complemented distributive lattice with a bottom element. A Boolean algebra with a top element will be called a {\em unital Boolean algebra}.
A homomorphism $\varphi: B_1\to B_2$ of Boolean algebras is called {\em proper} \cite{D}, provided that for any $b\in B_2$ there exists $a\in B_1$ such that $\varphi(a)\geq b$. In this paper any homomorphism of Boolean algebras is assumed to be proper. By $\GBA$ we denote the category of Boolean algebras and proper homomorphisms.  

By a {\em Boolean space} we mean what is usually called a locally compact Boolean space, that is a Hausdorff space in which compact-open sets form a base of the topology. Note that {\em compact Boolean spaces} are usually referred to as Boolean spaces. By $\LCBS$ we denote the category of Boolean spaces and continuous proper maps (recall that a map of topological spaces is {\em proper} if inverse images of compact sets are compact sets). 

Given a Boolean space $X$, all continuous maps $X\to\{0,1\}$, where $\{0,1\}$ is a discrete topological space, such that $f^{-1}(1)$ is a compact set, form a Boolean algebra denoted by $X^{*}$ and called the {\em dual Boolean algebra} of $X$.  The assignment $X\to X^*$ is the object part of the contravariant functor $\Al: \LCBS\to \GBA$. The restriction of $\Al$ to the category of compact Boolean spaces is an enriched contravariant $\Hom$-set functor, because given a compact Boolean space $X$,  the unital Boolean algebra $X^{*}$ consists of {\em all}\, continuous maps $X\to\{0,1\}$ (note that such a map is automatically also proper).

The {\em spectrum} $B^{*}$ of a Boolean algebra $B$ is the set of all non-zero homomorphisms from $B$ to the two-element Boolean algebra ${\bf 2}=\{0,1\}$. The set $B^{*}$ is equipped with the  topology whose base is formed by the sets $M(a)=\{f\in B^{*}: f(a)=1\}$, $a\in B$.  This space is a Boolean space and is called the {\em dual space} of $B$. It is well-known that the topology on $B^*$ can be also characterized as the subspace topology of the product space $\{0,1\}^B$. The assignment $B\to B^*$ is the object part of the contravariant functor $\Sp:\GBA\to\LCBS$. It is important for us that the restriction of $\Sp$ to the category of unital Boolean algebras is an enriched contravariant $\Hom$-set functor, because given a unital Boolean algebra $B$, the points of the space $B^{*}$ are {\em all}\,  unital Boolean algebra homomorphisms $B\to {\bf 2}$.

{\em Stone duality for unital Boolean algebras} \cite{S} (see also textbooks \cite{BS,GH}) states that the above described contravariant $\Hom$-set functors to $\{0,1\}$  establish a dual equivalence  between the categories of unital Boolean algebras and compact Boolean spaces. Therefore, $\{0,1\}$ is called a {\em dualizing object}, and this duality is {\em induced by a dualizing object}. 

{\em Stone duality for Boolean algebras} \cite{S, D} is an extension of the above duality. It states that the functors $\Al:\LCBS\to \GBA$ and $\Sp:\GBA\to\LCBS$ establish a dual equivalence between the categories $\LCBS$ and $\GBA$. 

In this paper, we will work  at the locally compact and non-unital level of generality. We will construct functors, whose restrictions to appropriate compact-unital subcategories are enriched $\Hom$-set functors, in the same way as in the commutative case outlined in this subsection.

\subsection{Skew Boolean algebras}\label{s22}  For a detailed introduction to the theory of skew Boolean algebras we refer the reader to \cite{BL,L2,L3}. A {\em skew lattice} is an algebra $(S;\wedge,\vee)$ of type $(2,2)$ such that the operations $\wedge$ and $\vee$ are associative, idempotent and satisfy the absorption identities $x\wedge(x\vee y)=x=x\vee(x\wedge y)$ and $(y\vee x)\wedge x=x=(y\wedge x)\vee x$. The {\em natural partial order} $\leq$ on a skew lattice $S$ is defined by $x\leq y$ if and only if $x\wedge y=y\wedge x=x$ or, equivalently, $x\vee y=y\vee x=y$. A skew lattice $S$ is {\em symmetric} if $x\wedge y=y\wedge x$ if and only if $x\vee y=y\vee x$. An element $0$ of $S$ is called a {\em zero} if $x\wedge 0=0\wedge x=0$ for all $x\in S$. $S$ is {\em Boolean} if it is symmetric, has a zero element and each principal subalgebra $\lceil x\rceil=\{y\in S:y\leq x\}=x\wedge S\wedge x$ forms a Boolean lattice.

Let $S$ be a Boolean skew lattice and $x,y\in S$. The {\em relative complement} $x\setminus y$ is the complement of $x\wedge y\wedge x$ in the Boolean lattice  $\lceil x\rceil$. A {\em skew Boolean algebra} is a Boolean skew lattice, whose signature is enriched by the nullary operation $0$ and the binary relative complement operation, that is, it is an algebra $(S;\wedge,\vee,\setminus,0)$. Skew Boolean algebras satisfy distributivity laws $x\wedge(y\vee z)=(x\wedge y)\vee(x\wedge z)$ and $(y\vee z)\wedge x=(y\wedge x)\vee (z\wedge x)$ \cite[2.5]{L6}.

Let $S$ be a skew lattice. It is called {\em rectangular} if there exist two sets $L$ and $R$  such that $S=L\times R$, and the operations $\wedge$ and $\vee$ are defined by $(a,b)\wedge (c,d)=(a,d)$ and $(a,b)\vee (c,d)=(c,b)$. Let $\D$ be the equivalence relation on $S$ given by $x\D y$ if and only if $x\wedge y\wedge x=x$ and $y\wedge x\wedge y=y$. It is known \cite[1.7]{L1} that $\D$ is a congruence relation, the $\D$-classes of $S$ are maximal rectangular subalgebras and the quotient $S/\D$ is the maximal lattice image of $S$. If $S$ is a skew Boolean algebra then $S/\D$ is the maximal Boolean algebra image of $S$ \cite[3.1]{BL}.

A skew lattice is called {\em left-handed} ({\em right-handed}) if it satisfies the identities $x\wedge y\wedge x=x\wedge y$ and $x\vee y\vee x=y\vee x$ (respectively, $x\wedge y\wedge x=y\wedge x$ and $x\vee y\vee x=x\vee y$). In a left-handed skew Boolean algebra the rectangular subalgebras are {\em flat} in the sense that $x\D y$ if and only if $x\wedge y=x$ and $y\wedge x=y$.
If $S$ is a left-handed skew Boolean algebra then the band $(S,\wedge)$ is left normal. It can be shown (or see \cite{KLa}) that for each $a\in S$ and $D\in S/{\mathcal D}$ such that $[a]_{\mathcal D}\geq D$ there is a unique $b\in S$ such that $[b]_{\mathcal D}=D$ and $a\geq b$. This element $b$ is called the {\em restriction} of $a$ to $D$ and is denoted by $a|_D$.

A skew Boolean algebra $S$ is called {\em primitive} if it has only one non-zero $\D$-class or, equivalently, if $S/\D$ is the Boolean algebra ${\bf 2}$. Up to isomorphism, finite primitive left-handed skew Boolean algebras are the  algebras ${\bf n+2}$, $n\geq 0$. These algebras  play an important role in this paper. The underlying set of ${\bf n+2}$ is  the set $\{0,1,\dots,n+1\}$ and its non-zero $\D$-class is $D=\{1,\dots,n+1\}$. The operations on $D$ are determined by lefthandedness: $x\wedge y=x$ and $x\vee y=y$ for any
$x,y\in D$. 

A skew Boolean algebra has {\em finite intersections} if any finite set of its elements has the greatest lower bound called the {\em intersection} with respect to the natural partial order. Since all intersections that we consider are finite, we will sometimes write just `intersections' for `finite intersections'. Homomorphisms of skew Boolean algebras with finite intersections are required to preserve the finite intersections. 

Let $\varphi: S_1\to S_2$ be a homomorphism of skew Boolean algebras and let $\overline{\varphi}:S_1/\D\to S_2/\D$ be the underlying homomorphism of Boolean algebras. We call $\varphi$ {\em proper} if $\overline{\varphi}$ is proper.

We fix the notation $\LSBA$ for the category of left-handed skew Boolean algebras and their are proper homomorphisms. All skew Boolean algebras, considered in the sequel, are left-handed, so we take a convention to write `skew Boolean algebra' for `left-handed skew Boolean algebra'. By a morphism of skew Boolean algebras we will mean a morphism in the category $\LSBA$.

\subsection{\'{E}tale spaces}\label{s23}

Preliminaries on \'{e}tale spaces can be found in any textbook on sheaf theory, e.g. in \cite{Br,MM}. An {\em \'{e}tale space over $X$} is a triple $(E,f,X)$, where $E$, $X$ are topological spaces and $f:E\to X$ is a surjective local homeomorphism.  The points of $E$ are called {\em germs.} A {\em local section} or just a {\em section} in $E$ is an open subset $A$ of $E$ such that the restriction of the map $f$ to $A$ is injective. If $U$ is an open set in $X$ then $E(U)$ is the set of all {\em sections} of $E$ over $U$, where a section $A$ is {\em over} $U$ provided that $f(A)=U$.  For $x\in X$ we denote the set of all $y\in E$ such that $f(y)=x$ by $E_x$. This set is called the {\em stalk} over $x$. If $A\in E(U)$ then for $x\in U$ by $A(x)$ we denote the germ $y\in A\cap E_x$.

Let $({\mathcal{A}},g,X)$ and $({\mathcal{B}},h,Y)$ be \'{e}tale spaces and $f:X\to Y$ be a continuous map. A {\em cohomomorphism over} $f$ (or an {\em $f$-cohomomorphism}) $k:{\mathcal{B}}\rightsquigarrow {\mathcal{A}}$ is a collection of maps $k_x: {\mathcal{B}}_{f(x)}\to {\mathcal{A}}_x$ for each $x\in X$ such that for every section $s\in {\mathcal{B}}(U)$ the function $x\mapsto k_x(s(f(x)))$ is a section of ${\mathcal{A}}$ over $f^{-1}(U)$. The maps $k_x$ are called the {\em components} of $k$.

We introduce the notation $\ESLCBS$ for the category of  \'{e}tale spaces over Boolean spaces and their cohomomorphisms over continuous proper maps.
All \'{e}tale spaces, considered in the sequel, are over Boolean spaces, therefore we
take a convention to write `\'{e}tale space' for `\'{e}tale space over a Boolean space'. By a morphism of \'{e}tale spaces we will mean a morphism in the category $\ESLCBS$.

\subsection{Equivalence of the categories of skew Boolean algebras and \'{e}tale spaces.}\label{s24}
The constructions of this paper rely on the generalization of Stone duality to left-handed skew Boolean algebras established in \cite{K}, that is outlined in this subsection.  For a detailed exposition and proofs we refer the reader to \cite{K}.

 Let $(E,\pi, X)$ be an \'{e}tale space. Denote by $E^{*}$ the set of all compact-open sections of $E$. Let $A,B\in E^*$.  We put $$(A\ucup B)(x)=\left\lbrace\begin{array}{ll}B(x), & \text{ if } x\in \pi(B),\\
A(x), & \text{ if } x\in \pi(A)\setminus \pi(B),
\end{array}\right.$$

$$ (A\ocap B)(x)=A(x) \text{ for all } x\in \pi(A)\cap \pi(B).$$

Then $\ucup$ and $\ocap$ are well-defined binary operations on $E^*$ and
$(E^*, \ucup, \ocap)$ is a left-handed Boolean skew lattice.  
The associated skew Boolean algebra $E^{*}$ is called the {\em dual skew Boolean algebra} to the \'{e}tale space $(E,\pi,X)$. Note that $a{\mathcal D}b$ in $E^*$ if and only if $\pi(a)=\pi(b)$ and the assignment $[a]_{\mathcal D}\mapsto \pi(a)$ establishes an isomorphism between $(S/{}\mathcal D)^*$ and $X^*$.
Remark that $E^*$ has finite intersections if and only if $E$ is Hausdorff.

In the opposite direction, let $S$ be a skew Boolean algebra. The points of the dual \'{e}tale space of $S$ are the ultrafilters of $S$. The latter can be characterized as follows. Let $F$ be an ultrafulter of the Boolean algebra $S/{\mathcal D}$ and let $a\in S$ be such that $[a]_{\mathcal D}\in F$. Then the set
$$X_{a,F}=\{b\in S\colon \text{ there is } c\in S \text{ such that }a,b\geq c \text{ and } [c]_{\mathcal D}\in F\}
$$
is an ultrafilter in $S$, and any ultrafilter in $S$ is of this form. 

The set  $S^{*}$  of all ultrafilters of $S$ is called the {\em spectrum} of $S$. Let $\hat{\pi}:S^{*}\to(S/\D)^{*}$ be the {\em projection  map} given by
$X_{a,F} \mapsto F.$

For $a\in S$ we put
$$
M(a)=\{F\in S^{*}:a\in F\}.
$$
The topology on $S^{*}$ is given by a subbase,  that is in fact a base, consisting of the sets $M(a)$, $a\in S$.  We have that $S^{*}=(S^{*},\hat{\pi},(S/\D)^{*})$ is an \'{e}tale space called the {\em dual \'{e}tale space} of $S$. 
The space $S^*$ is Hausdorff if and only if the skew Boolean algebra $S$ has finite intersections.

We now describe the correspondence for morphisms. 
Let $(E,f,X),(G,g,Y)$ be \'{e}tale spaces and $k:E\to F$ a morphism.  We extend $k$  to the map on sections. This map takes compact-open sections to compact-open sections and preserves the operations $\ucup$, $\ocap$ and $\varnothing$.  This leads to the morphism $\widetilde{k}:E^*\to G^*$ of skew Boolean algebras. 
 Note that if the spaces $E$ and $G$ are Hausdorff and all components of $k$ are injective then
$\widetilde{k}$ preserves finite intersections.

Conversely, let $S,T$ be skew Boolean algebras and  $k:S\to T$ a morphism. The map $\overline{k}^{-1}$ induces a morphism of Boolean spaces from $(T/\D)^{*}$ to $(S/\D)^{*}$. 
Let $F\in (T/\D)^{*}$ and $V\in T^{*}_F$. Then the set $k^{-1}(V)$ is either empty or a disjoint union of several ultrafilters over ${\overline{k}^{-1}(F)}$. This allows to define the component $\widetilde{k}_F$ of the morphism $\widetilde{k}:S^*\to T^*$. The domain of $\widetilde{k}_F$ is $S^*_{\overline{k}^{-1}(F)}$ and for $U\in S^{*}_{\overline{k}^{-1}(F)}$ we set $\widetilde{k}_F(U)=V$ if  $V\in T^{*}_F$ and $U\subseteq k^{-1}(V)$. Note that if $S,T$ have finite intersections and $k$ preserves them then all components of $\widetilde{k}:S^*\to T^*$ are injective.

We now state the theorem that provides insights needed for establishing the main results of this paper.

\begin{theorem}[\cite{K}]\label{th1} The described assignments are functors that establish an equivalence between the categories $\ESLCBS$ and $\LSBA$.  If $S$ is a skew Boolean algebra then $S$ is naturally isomorphic to $S^{**}$ via the map $\beta_S$ given by
$\beta_S(a)=M(a)$, $a\in S$. If $E$ is an \'{e}tale space over a Boolean space then it is naturally isomorphic to $E^{**}$ via the map $\gamma_E$ given by $
\gamma_E(A)=N_A=\{N\in E^{*}:A\in N\}$, $A\in E$. \end{theorem}

Consider an application of this theorem. Let $n\geq 0$ and $\varphi:S\to {\bf n+2}$ be a morphism of skew Boolean algebras (such morphisms will play an important role in Section \ref{s6}). We interpret this topologically.  Observe that $({\bf n+2})/{\mathcal D}={\bf 2}$ and ${\bf 2}^{*}$ is a one-element space, $\{a\}$. The cohomomorphism $\widetilde{\varphi}$ has therefore only one component $\widetilde{\varphi}_a: S^*_{\overline{\varphi}^{-1}(a)}\to \{1,\dots, n+1\}$, that completely determines $\varphi$. Moreover, any map $\phi: S^*_{\overline{\varphi}^{-1}(a)}\to \{1,\dots, n+1\}$ gives rise to a morphims from $S$ to ${\bf n+2}$ over $\overline{\varphi}^{-1}$.

We also have the following theorem.

\begin{theorem}[\cite{K, BCV}] \label{th:a}The category of skew Boolean algebras with intersections is equivalent to the category of Hausdorff \'{e}tale spaces over Boolean spaces whose morphisms are cohomomorphisms with injective components.
\end{theorem}

 Let $S,T$ be skew Boolean algebras with intersections. Theorems \ref{th1} and \ref{th:a} tell us that, unless $S$ is commutative, the set of all morphisms from $S$ to $T$ is much bigger than the set of intersection-preserving such morphisms.  For example, if $S={\bf 4}$ and $T={\bf 3}$ morphisms from $S$ to $T$ are given by functions
from $\{1,2,3\}$ to $\{1,2\}$.  Intersection-preserving morphisms are given by injective such functions, which means that there are no intersection-preserving morphisms at all (note that we do not take into account the zero map since we consider only proper morphisms). For the constructions of this paper, it is crucial that we consider all cohomomorphisms between appropriate spaces, and not only ones with all components injective. So it is important that the duality theorem we use in this paper is Theorem \ref{th1} and not Theorem \ref{th:a}.

We now explain why the unital versions of dualities of Theorems \ref{th1} and \ref{th:a} are not induced by a dualizing object. For Theorem \ref{th1}, consider the functor constructing the spectrum of a skew Boolean algebra. As follows from Lemma 6.3 of \cite{K} points of the spectrum of a skew Boolean algebra $S$ are in a bijective correspondence with morphisms $S\to {\bf n+2}$ such that the inverse image of $1$ is non-empty and minimal possible. This does not produce all morphisms $S\to {\bf n+2}$ (and neither all morphisms from $S$ to another skew Boolean algebra), even in the case when $S/{\mathcal D}$ is unital. For the example from the previous paragraph, there are exactly $2^3=8$ morphisms from ${\bf 4}$ to ${\bf 3}$. But only three of these morphisms, that are listed below, give rise to the points of the spectrum of ${\bf 4}$:
$$
\varphi_1: \begin{array}{l}0\mapsto 0\\1\mapsto 1\\
 2 \mapsto 2\\
 3\mapsto 2\end{array}; \,\,\,\,\,
 \varphi_2:\begin{array}{l}0\mapsto 0\\1\mapsto 2\\
 2 \mapsto 1\\
 3\mapsto 2\end{array}; \,\,\,\,\,
 \varphi_3:\begin{array}{l}0\mapsto 0\\1\mapsto 2\\
 2 \mapsto 2\\
 3\mapsto 1\end{array}.$$
 For Theorem \ref{th:a}, this follows from the fact that for any skew Boolean algebra with intersections $S$, there exists another skew Boolean algebra with intersections $T$ such that there are no intersection-preserving morphisms from $S$ to $T$ at all (we have only to take care that any stalk of $T^*$ has cardinality sticktly less than the cardinality of any stalk of $S^*$).

\section{The functors $\lambda_n$, $n\geq0$}\label{s4}

Let $X$ be a Boolean space and let $n\geq 0$ be fixed. We regard the set $\{0,\dots,n+1\}$ as a discrete topological space. Let $\lambda_n(X)$ denote the set of all continuous maps $f$ from $X$ to $\{0,\dots,n+1\}$ such that the sets$f^{-1}(1)$, $\dots$, $f^{-1}(n+1)$ are compact. Define the binary operations $\wedge$ and $\vee$ on $\lambda_n(X)$ to be induced by the operations  $\wedge$ and $\vee$ on the primitive skew Boolean algebra ${\bf n+2}$. That is, for $f,g\in\lambda_n(X)$ we put
$$
(f\wedge g)(x)=f(x)\wedge g(x), \,\, (f\vee g)(x)=f(x)\vee g(x).
$$

With respect to $\wedge$ and $\vee$ the set $\lambda_n(X)$ becomes a left-handed Boolean skew lattice. By adding to its signature the relative complement operation and the zero, we turn itno into a left-handed skew Boolean algebra. Note that  $\lambda_0(X)=X^*$. 

It is immediate that the elements of $\lambda_n(X)$ are in a bijective correspondence with ordered $(n+1)$-tuples of pairwise disjoint compact-open subsets of of $S$ via the assignment
$$
f\mapsto (f^{-1}(1),\dots, f^{-1}(n+1)).
$$

\begin{remark}\label{rem:ls1} Yet another realization of $\lambda_n(X)$ is by flags $A_{n+1}\supseteq \dots \supseteq A_1$ of compact-open subsets of $X$ via the assignment 
$$
f\mapsto f^{-1}(\{1,\dots,n+1\})\supseteq f^{-1}(\{1,\dots,n\})\supseteq f^{-1}(1).
$$
Let $A_{n+1}\supseteq \dots \supseteq A_1$ and  $B_{n+1}\supseteq \dots \supseteq B_1$ be the flags corresponding to $f$ and $g$, respectively. It is easy to verify, applying  the definitions of $\vee$ and $\wedge$, that the flag, corresponding to $f\vee g$ is $C_{n+1}\supseteq \dots \supseteq C_1$, where $C_i=(A_i\setminus B_{n+1})\cup B_{i}$ for all $i$, 
and the flag, corresponding to $f\wedge g$ is $D_{n+1}\supseteq \dots \supseteq D_1$, where $D_i=A_i\cap B_{n+1}$ for all $i$. From this description it follows that  $\lambda_1(X)=\omega(\Al(X))$, where $\Al$ is the functor from Subsection \ref{s21} and $\omega: \GBA\to\LSBA$ is the (version with proper morphisms of the) Leech-Spinks $\omega$-functor \cite{LS}. Therefore, the construction of $\lambda_1$ provides a natural interpretation for the  functor $\omega$. 
\end{remark}

Let us look at the structure of the skew Boolean algebra $\lambda_n(X)$ in more detail. For $f\in \lambda_n(X)$ by ${\hat{f}}: X\to\{0,1\}$ we denote the map given by $${\hat{f}}^{-1}(1)=f^{-1}(\{1,\dots, n+1\}).$$

\begin{lemma} Let $f,g\in \lambda_n(X)$. Then $f{\mathcal D} g$ if and only if $f^{-1}(\{1,\dots, n+1\})=g^{-1}(\{1,\dots, n+1\})$. The Boolean algebra $\lambda_n(X)/{\mathcal{D}}$ is isomorphic to the Boolean algebra $\lambda_0(X)=X^*$ via the map $[f]_{\mathcal D}\mapsto {\hat{f}}$. 
\end{lemma}

\begin{proof} By definition,  $f{\mathcal D} g$ is equivalent to $f\wedge g=f$ and $g\wedge f=g$. Therefore, $f(x)\wedge g(x)=f(x)$ and $g(x)\wedge f(x)=g(x)$ for all $x\in X$. Applying the description of the operation $\wedge$ on ${\bf n+2}$ we see that the latter equalities are equivalent to the condition $f(x)=0$ if and only if $g(x)=0$ for all $x\in X$. This is clearly equivalent to  $f^{-1}(\{1,\dots, n+1\})=g^{-1}(\{1,\dots, n+1\})$. The second statement is straightforward to verify.
\end{proof}

\begin{lemma} \label{lem:order} The natural partial order in $\lambda_n(X)$ is given by $f\geq g$ if and only if $f^{-1}(i)\supseteq g^{-1}(i)$ for all $i=1,\dots, n+1$.
\end{lemma}

\begin{proof} By definition, $f\geq g$ is equivalent to $f\wedge g=g\wedge f=g$.
By the definition of the operation $\wedge$ on $\lambda_n(X)$, the latter is equivalent to $f(x)\wedge g(x)=g(x)\wedge f(x)=g(x)$ for all $x\in X$, where the equalities are in ${\bf n+2}$. This is equivalent to the condition $g(x)=i$ implies $f(x)=i$ for all $i=1,\dots, n+1$.
\end{proof}

\begin{lemma} The skew Boolean algebra $\lambda_n(X)$ has finite intersections. If $f,g\in \lambda_n(X)$ then their intersection $f\cap g$ is given by $(f\cap g)^{-1}(i)=f^{-1}(i)\cap g^{-1}(i)$ for all $i=1,\dots, n+1$.
\end{lemma}

\begin{proof} Let $h\in\lambda_n(X)$ is given by  $h^{-1}(i)=f^{-1}(i)\cap g^{-1}(i)$ for all $i=1,\dots, n+1$. Then $f,g\geq h$ by Lemma \ref{lem:order}. Assume that $f,g\geq q$ and let $i\in\{1,\dots, n+1\}$. Then $f^{-1}(i)\supseteq q^{-1}(i)$ and $g^{-1}(i)\supseteq q^{-1}(i)$. It follows that $h^{-1}(i)\geq q^{-1}(i)$. This proves that $f\cap g$ exists and equals $h$.
\end{proof}

It is well-known that ultrafilters of the Boolean algebra $\lambda_0(X)$ are the sets
$\{f\in \lambda_0(X)\colon f(x)=1\}$. Consequently ultrafilters of the Boolean algebra $\lambda_n(X)/{\mathcal{D}}$ are the sets
\begin{equation}\label{eq:nx}
N_x=\{[f]_{\mathcal D}\in \lambda_n(X)/{\mathcal{D}}\colon f(x)\in \{1,\dots, n+1\}\}, \,\, x\in X.
\end{equation}

Let $\alpha: \lambda_n(X)\to  \lambda_n(X)/{\mathcal{D}}$ be the projection map.

\begin{lemma} \label{lem:jun15} For every $x\in X$ and $i\in\{1,\dots, n+1\}$ the set
$$
N_{x,i}=\{f\in \lambda_n(X)\colon f(x)=i\}
$$
is an ultrafilter of $\lambda_n(X)$, and any ultrafilter is of this form. In addition, $\alpha(N_{x,i})=N_x$.
\end{lemma}

\begin{proof} Let $f\in N_{x,i}$.  Show that $N_{x,i}=X_{f,N_x}$.
Let $g\in N_{x,i}$. Since $x\in f^{-1}(i)\cap g^{-1}(i)$ then $f\cap g\in N_{x,i}$.
Thus also $[f\cap g]_{\mathcal D}\in N_x$. This and $f,g\geq f\cap g$ imply that $g\in X_{f,N_x}$. 

Conversly, let $g\in X_{f,N_x}$. Then there is some $h\in\lambda_n(X)$ with $[h]_{\mathcal D}\in N_x$ such that $f,g\geq h$. Since $f\geq h$ and $h(x)\neq 0$ it follows that $h(x)=i$ by Lemma \ref{lem:order}. Now, $g\geq h$ implies $g(x)=i$ again by Lemma \ref{lem:order}. It follows that $g\in N_{x,i}$, as required. The second claim is immediate.
\end{proof}

\begin{corollary} For any ultrafilter $N_x$ of  $\lambda_n(X)/{\mathcal{D}}$ we have
$$
\alpha^{-1}(N_x)=\bigcup_{i=1}^{n+1}N_{x,i}.
$$
Consequently, each stalk of the \'{e}tale space $(\lambda_n(X))^*$ has cardinality $n+1$.
\end{corollary}
\begin{proof} Ultrafilters in $\lambda_n(X)$ that are over $N_x$ are of the form
$X_{f,N_x}$, where $f(x)\in \{1,\dots, n+1\}$. The statement now follows from the equality $N_{x,i}=X_{f,N_x}$, where $f(x)=i$, established in the proof of Lemma \ref{lem:jun15}.
\end{proof}

Let $X^{(i)}=\{N_{x,i}\}_{x\in X}$, $i\in\{1,\dots, n+1\}$. It follows from Lemma \ref{lem:jun15} that, as a set, $(\lambda_n(X))^{*}$ is the union
$X^{(1)}\cup\dots\cup X^{(n+1)}$ of $(n+1)$ disjoint copies of $X$. 

\begin{lemma} The space $(\lambda_n(X))^{*}$ as a topological space is a disjoint union of $(n+1)$ copies of the Boolean space $X$.
\end{lemma}

\begin{proof}   
Consider each $X^{(i)}$ as a topological space homeomorphic to $X$ via the map
$N_{x,i}\mapsto x$. Then $(\lambda_n(X))^{*}$ can be considered as a disjoint union of the spaces $X^{(i)}$. We aim to show that this topology on $(\lambda_n(X))^{*}$ coincides with the dual \'{e}tale space topology, whose construction is outlined in Subsection \ref{s24}. To show this, it is enough to verify that the map $\pi:(\lambda_n(X))^{*} \to X$ given by $N_{x,i}\mapsto x$ is a local homeomorphism. Fix some $N_{x,i}$ and a compact-open set $A$ in $X$ such that $x\in A$. Consider the function $f\in\lambda_n(X)$ given by $f^{-1}(i)=A$ and $f^{-1}(j)=\varnothing$, $j\in\{1,\dots, n+1\}\setminus\{i\}$. It is clear that $M(f)=\cup_{y\in A}N_{y,i}$. This is a neighborhood of $N_{x,i}$ that is homeomorphic via $\pi$ to $A$ since a basic open subset of $M(f)$ equals $M(g)$, $f\geq g$, and $\pi(M(g))=B\subseteq A$, where $g^{-1}(i)=B$. \end{proof}

We now define $\lambda_n$ on morphisms.  Let $n\geq 0$, $g: X_1\to X_2$ be a morphism of Boolean spaces and
$f\in \lambda_n(X_2)$. We put
\begin{equation}\label{eq6}
\lambda_n(g)(f)=fg
\end{equation}

It is easy to check that $\lambda_n(g)$ is a skew Boolean algebra morphism from $\lambda_n(X_2)$ to $\lambda_n(X_1)$ and that this makes $\lambda_n$ a contravariant functor from the category $\LCBS$ to the category $\LSBA$.

We finish this section by recording the following fact that will be needed in Section \ref{s7}.

\begin{lemma}\label{lem:aaa} Let $g: X_1\to X_2$ be a morphism of Boolean spaces. Then 
$\overline{\lambda_n(g)}^{-1}(N_x)=N_{g(x)}$ 
and $\lambda_n(g)^{-1}(N_{x,i})=N_{g(x),i}$
for any $x\in X_1$ and $i=1,\dots, n+1$.\end{lemma}

\begin{proof} We have $\overline{\lambda_n(g)}([f]_{\mathcal D})=[fg]_{\mathcal D}$. Therefore, $\overline{\lambda_n(g)}([f]_{\mathcal D})\in N_x$ if and only if
$[fg]_{\mathcal D}\in N_x$. By \eqref{eq:nx} this is equivalent to $fg(x)\in \{1,\dots, n+1\}$. The latter means that $[f]_{\mathcal D}\in N_{f(x)}$ and implies the first equality. For the second equality, observe that
$$
\lambda_n(g)(f)\in N_{x,i} \,\Leftrightarrow \, fg(x)\in N_{x,i} \,\Leftrightarrow \, f\in N_{g(x),i}.
$$
\end{proof}

\section{The functors $\Lambda_n$ and the adjunctions $\Lambda_n \dashv \lambda_n$, $n\geq 0$}\label{s6}

Let $S$ be a skew Boolean algebra, $\alpha:S\to S/\D$ the canonical projection and $\hat{\alpha}$ the projection from $S^{*}$ to $(S/\D)^{*}$ given by $\hat{\alpha}(X_{a,F})=F$. We also fix $n\geq 0$. 

Let $\{0,1,\dots, n+1\}^{S}$ be the set of all maps from $S$ to  $\{0,1,\dots,n+1\}$. We consider $\{0,1,\dots, n+1\}$ as a discrete space and $\{0,1,\dots, n+1\}^{S}$ as a product space. That is, a base of the topology on $\{0,1,\dots,n+1\}^{S}$ is formed by the sets
\begin{equation}\label{jun12}
U_{\delta}=\bigcap_{t\in T}\{f\in \{0,1,\dots,n+1\}^{S}: f(t)=\delta(t)\},
\end{equation}
where $T$ runs through the finite subsets of $S$ and $\delta$ runs through the functions from $T$ into $\{0,1,\dots,n+1\}$.

We denote by $\Lambda_n(S)$ the set of all morphisms from $S$ to ${\bf n+2}$ in the category $\LSBA$. We endove $\Lambda_n(S)$ with the subspace topology inherited from the product topology on the space $\{0,1,\dots,n+1\}^{S}$.

For $s\in S$ and $i\in \{0,1,\dots, n+1\}$  put 
$$
L(s,i)=\{f\in \Lambda_n(S)\colon f(s)=i\}.
$$

\begin{lemma} \label{lem:jun13a} The sets $L(s,i)$, where $s$ runs through $S$ and $i\in\{1,\dots, n+1\}$, form a subbase of the topology on $\Lambda_n(S)$.
\end{lemma}

\begin{proof} By the definition of subspace topology and in view of the base given by \eqref{jun12} we obtain that the sets
$$
V_{\delta}=U_{\delta}\cap \Lambda_n(S) = \bigcap_{t\in T}\{f\in \Lambda_n(S): f(t)=\delta(t)\},
$$
where $T$ runs through the finite subsets of $S$ and $\delta$ runs through the functions from $T$ into $\{0,1,\dots,n+1\}$, form a base of the topology on $\Lambda_n(S)$. Therefore, the topology on $\Lambda_n(S)$ admits a subbase consisting from the sets $L(s,i)$, $s\in S$, $i\in \{0,1,\dots, n+1\}$.  Let $s\in S$. Show that the set $L(s,0)$ can be expressed as a union of some of the sets $L(t,i)$, where $t\in S$ and $i\in \{1,\dots, n+1\}$. Let $A=\{t\in \Lambda_n(S)\colon  \alpha(t)\wedge \alpha(s)=0\}$. We aim to show that
\begin{equation}\label{eq:jun13}
L(s,0)=\bigcup_{t\in A}\bigcup_{i=1}^{n+1}L(t,i).
\end{equation}

Let $f\in L(s,0)$. Since $f$ is non-zero, there is $t\in S$ such that $f(t)\neq 0$.
Observe that $f(s|_{\alpha(s)\wedge \alpha(t)})=0$, since $s\geq s|_{\alpha(s)\wedge \alpha(t)}$ and $f(s)=0$. It follows that $f(t|_{\alpha(s)\wedge \alpha(t)})=0$ as well since the two elements $s|_{\alpha(s)\wedge \alpha(t)}$ and $t|_{\alpha(s)\wedge \alpha(t)}$ are in the same ${\mathcal D}$-class. It follows that
$f(t|_{\alpha(t)\setminus \alpha(s)})\neq 0$. We obtain that $t|_{\alpha(t)\setminus \alpha(s)}\in A$ and 
$f\in \bigcup_{i=1}^{n+1}L(t|_{\alpha(t)\setminus \alpha(s)},i)$, proving that
$L(s,0)\subseteq\bigcup_{t\in A}\bigcup_{i=1}^{n+1}L(t,i)$.

To prove the reverse inclusion, we let $f\in \bigcup_{t\in A}\bigcup_{i=1}^{n+1}L(t,i)$.  We have $\alpha(s\wedge t)=\alpha(s)\wedge \alpha(t)=0$, so that $s\wedge t=0$. It follows that $f(s\wedge t)=f(0)=0$. On the other hand we have $f(s\wedge t)=f(s)\wedge f(t)$. We obtain that $0=f(s)\wedge f(t)$ and $f(t)\neq 0$. It follows that $f(s)=0$ and thus $f\in L(s,0)$.  This finishes the proof of \eqref{eq:jun13}.
\end{proof}

\begin{lemma}\label{lem:jun13} Let $s\in S$ and $i\in\{1,\dots, n+1\}$. The set $L(s,i)$ is a closed subset of the space $\{0,\dots,n+1\}^S$. Consequently, $L(s,i)$ is a compact-open subset of $\Lambda_n(S)$.
\end{lemma}
\begin{proof} Let $Y_{\wedge}$ and $Y_{\vee}$ denote the sets of functions from $S$ to $\{0,1,\dots, n+1\}$ that preserve $\wedge$ or $\vee$, respectively.
Let, further, $Y_0$ and $Y_{s,i}$ denote the sets  of functions from $S$ to $\{0,1,\dots, n+1\}$ that map $0$ to $0$ or that map $s$ to $i$, respectively.
It can be shown applying a standard argument (see, for example, proof of Lemma 1, Chapter 34, p. 326 of \cite{GH}) that each of the sets $Y_{\wedge}$, $Y_{\vee}$, $Y_0$ and $Y_{s,i}$ is a  closed subset of $\{0,1,\dots, n+1\}^S$. It follows that $L(s,i)$ is closed, too, since $L(s,i)=Y_{\wedge}\cap Y_{\vee}\cap Y_{0}\cap Y_{s,i}$. 

For the second claim, observe that $L(s,i)$ is compact in $\{0,\dots,n+1\}^S$ as a closed subset of a compact space. Therefore, $L(s,i)$ is compact-open in $\Lambda_n(S)$.\end{proof}

\begin{theorem} $\Lambda_n(S)$ is a Boolean space.
\end{theorem}
\begin{proof}  $\Lambda_n(S)$ is Hausdorff as a subspace of the Hausdoeff space $\{1,\dots, n+1\}^S$. So, in view of Lemmas \ref{lem:jun13a} and \ref{lem:jun13}, we have that $\Lambda_n(S)$ is a  Hausdorff space, in which  compact-open sets form a base of the topology. Thus $\Lambda_n(S)$ is a Boolean space.
\end{proof}

We now define $\Lambda_n$ on morphisms.  Let $h:S_1\to S_2$ be a morphism of skew Boolean algebras. For $f\in \Lambda_n(S_2)$ we set
\begin{equation}\label{aux8}
(\Lambda_n(h))(f)=fh.
\end{equation}

\begin{lemma}\label{lem24} $\Lambda_n(h)$ is a morphism of Boolean spaces from  $\Lambda_n(S_2)$ to $\Lambda_n(S_1)$.
\end{lemma}

\begin{proof} It is immediate that for each $f\in \Lambda_n(S_2)$ we have that
 $(\Lambda_n(h))(f)\in \Lambda_n(S_1)$. To show that $\Lambda_n(h)$ is a continuous proper map,  it is enough to verify that the set $(\Lambda_n(h))^{-1}(L(s,i))$ is  compact-open for any $s\in S_2$ and any $i\in\{1,\dots, n+1\}$. We observe that
\begin{multline*}
f\in (\Lambda_n(h))^{-1}(L(s,i)) \Longleftrightarrow  (\Lambda_n(h))(f)\in L(s,i) \Longleftrightarrow \\
fh\in L(s,i) \Longleftrightarrow f(h(s))=i \Longleftrightarrow f\in L(h(s),i),
\end{multline*}
 implying that  $(\Lambda_n(h))^{-1}(L(s,i))=L(h(s),i)$, and the statement follows.\end{proof}

It is straightforward to check that the constructed assignments define a contravariant functor $\Lambda_n:\LSBA\to\LCBS$.

We are now prepared to formulate and prove our adjunction theorem.

\begin{theorem}\label{th25} For each $n\geq 0$ the functor $\Lambda_n:\LSBA\to\LCBS^{op}$ is the left adjoint to the functor $\lambda_n:\LCBS^{op}\to \LSBA$. For a skew Boolean algebra $S$ the component $\eta_S$ of the unit of the adjunction $\eta$ is given by 
\begin{equation}\label{aux13}\eta_S(a)(g)=g(a),  a\in S, g\in \Lambda_n(S).
\end{equation}
\end{theorem}

\begin{proof} It is immediate that $\eta_S$ is a morphism. Let  $S$ be a skew Boolean algebra, $X$ a Boolean space and  $\mu: S\to \lambda_n(X)$ a morphism of skew Boolean algebras. Our aim is to show that there is a unique morphism $u: X\to \Lambda_n(S)$ of Boolean spaces such that $\mu=\lambda_n(u)\eta_S$.

For each $x\in X$ we put
\begin{equation}\label{eq:jun20}
u(x)(s)=\mu(s)(x), s\in S.
\end{equation}

Let us verify that $u$ is a proper continuous map. For this, we have to show that the inverse image under $u$ of a compact-open subset of $\Lambda_n(S)$ is  compact-open in $X$. Since any compact-open subset is a finite union of basic compact-open sets, and any basic compact-open set is a finite intersection of the sets of the form $L(s,i)$,  it is enough to verify that $u^{-1}(L(s,i))$ is compact-open in $X$ for each $s\in S$ and each $i\in\{1,\dots,n+1\}$.
We have
\begin{multline*}
x\in u^{-1}(L(s,i)) \Leftrightarrow u(x)\in L(s,i) \Leftrightarrow u(x)(s)=i \Leftrightarrow \\
\mu(s)(x)=i \Leftrightarrow x\in \mu(s)^{-1}(i).
\end{multline*}
Since $\mu(s)\in \lambda_n(X)$ then $\mu(s)^{-1}(i)$ is compact-open by the definition of $\lambda_n(X)$. So $u^{-1}(L(s,i))$ is compact-open, too.

Verify the equality $\mu=\lambda_n(u)\eta_S$. For each $a\in S$ and $x\in X$ we have
$$\begin{array}{lclr}
(\lambda_n(u)\eta_S(a))(x)&=&\eta_S(a)u(x)&\text{(by \eqref{eq6})}\,\,\\
&=&u(x)(a)&\text{(by \eqref{aux13})}\,\,\\
&=&\mu(a)(x)&\text{(by \eqref{eq:jun20})},
\end{array}$$
as required.

The uniqueness of  $u$ is shown in a standard way, we leave the details to the reader.
\end{proof}

\begin{remark} \label{rem:ls2} Let $S$ be a skew Boolean algebra and $\Al$, $\Sp$ functors from Subsection \ref{s21}. We set $\Omega=\Al\Lambda_1:\LSBA\to \GBA$. It follows from theorem \ref{th25} that $\Omega$ is the left adjoint to the functor  $\omega=\lambda_1\Sp: \GBA\to \LSBA$ from \cite{LS}. The construction of $\Lambda_1$ yields the full description of (the version with proper morphisms of) $\Omega$, that had been questioned in \cite{LS}.
\end{remark}

The statement, that follows,  provides a way  to faithfully represent a skew Boolean algebra $S$ as a subalgebra of a well-understood skew Boolean algebra with intersections $\lambda_n\Lambda_n(S)$ whose dual space has $(n+1)$-element stalks. If $n=1$ we obtain a faithful representation of $S$ in a skew Boolean algebra $\lambda_1\Lambda_1(S)$ with intersections whose dual space has has $2$-element stalks. This algebra has, roughly speaking, a `low degree of non-commutativity' (but note that its underlying Boolean algebra $(\Lambda_1(S))^*$ is rather huge in comparison with the underlying Boolean algebra $S/{\mathcal D}$ of $S$).  
In the terminology of Leech and Spinks \cite{LS}, $\lambda_1\Lambda_1(S)$ is a {\em minimal skew Boolean cover} of the Boolean algebra $(\Lambda_1(S))^{*}$. We refer the reader to \cite{LS} for an interesting discussion of different approaches to the definition of the notion of a minimal skew Boolean cover.  

\begin{theorem}\label{th:3}
If $n\geq 1$ then the map $\eta_S$ is injective, and is therefore a faithful representation of
$S$ as a subalgebra of $\lambda_n\Lambda_n(S)$.
\end{theorem}

\begin{proof} Let $a,b\in S$ and $a\neq b$. This implies that $M(a)\neq M(b)$. Then we can assume that there is an ultrafilter $X_{a,F}$ in $S$ such that $b\not\in X_{a,F}$ (note that $a\in X_{a,F}$ by definition of $X_{a,F}$).  Let $x$ be the only point of the space ${\bf 2}^{*}$. The map $g:{\bf 2}^{*}\to (S/{\mathcal D})^*$ given by $x\mapsto F$ is a continuous map. As follows from the discussion in the paragraph after Theorem \ref{th1}, we can construct the cohomomorphism $k:S^*\to  ({\bf n+2})^*$ over $g$ given by $k_x(X_{a,F})=1$ and $k_x(G)=2$ for any $G \in S^*_F$, $G\neq X_{a,F}$. Since $b\not\in X_{a,F}$ then $k(b)\neq 1$. It follows that the morphism from $S$ to ${\bf n+2}$, that corresponds to $k$, has different values at $a$ and $b$. This implies that $\eta_S(a)(k)\neq \eta_S(b)(k)$.
\end{proof}

\section{The algebras of the monads of the adjunctions $\Lambda_n\dashv \lambda_n$}\label{s7}

Let $n\geq 0$ be fixed throughout this section. Preliminaries on monads can be found in, e.g., \cite[Chapter 10]{Awo} or \cite[Chapter VI]{M}. 

Let $(T,\eta,\mu)$ be the monad over the category $\LSBA$ that arises from the adjunction $\Lambda_n\dashv\lambda_n$. We have $T=\lambda_n\Lambda_n$, $\eta$ is the unit of the adjunction given by \eqref{aux13}, and $\mu=\lambda_n\epsilon_{\Lambda_n}$ is a natural transformation from $T^2$ to $T$, where $\epsilon: 1_{\LCBS}\to \Lambda_n\lambda_n$ is the counit of the adjunction.

The following statement is straightforward to verify.

\begin{lemma} \label{lem:counit} Let $X$ be a Boolean space.  Then 
\begin{equation}\label{eq:9jul}
\epsilon_X(x)(f)=f(x), x\in X, f\in \lambda_n(X).
\end{equation}
\end{lemma}

Let $S$ be a skew Boolean algebra. In this section we will need to work with skew  Boolean algebras $\lambda_n\Lambda_n\lambda_n\Lambda_n(S)$, $\lambda_n\Lambda_n(S)$ and morphisms between them. These objects are somewhat complicated, and to proceed, we first propose a convenient way of working with them. Our insights will be based on the topological representation of morphisms of skew Boolean algebras by morphisms of their dual spaces.
 
We start from a useful description of the points of the space $\Lambda_n(S)$.

\begin{lemma}\label{lem20}
There is a bijective correspondence between  the points of the space $\Lambda_n(S)$ and  elements of the set 
\begin{equation}\label{eq300}\{(F,f)\colon F\in (S/{\mathcal D})^{*}, f\in\{1,\dots,n+1\}^{S^{*}_F}\}.
\end{equation}
\end{lemma}

\begin{proof} Let $a$ be the only point of the Boolean space $(({\bf{n+2}})/{\mathcal D})^* = {\bf 2}^*$. Let $h\in\Lambda_n(S)$. Then $\overline{h}^{-1}$ induces a continuous proper map, that we denote also by $\overline{h}^{-1}$,  from  ${\bf 2}^*$ to $(S/{\mathcal D})^*$.  Let $F_h\in(S/\D)^{*}$ be such that $\overline{h}^{-1}(a)=F_h$. Then $\widetilde{h}$ has the only one component $\widetilde{h}_a: S^{*}_{F_h}\to {\bf{(n+2)}}^{*}_a=\{1,\dots,n+1\}$. Show that the constructed correspondence $h\mapsto (F_h, \widetilde{h}_a)$ is bijective. Assume that $(F_h, \widetilde{h}_a)= (F_g, \widetilde{g}_a)$. This immediately implies that  $\overline{h}=\overline{g}$ and that $h^{-1}(i)=g^{-1}(i)$ for all $i=1,\dots, n+1$. This proves that the constructed correspondence is injective. For the reverse direction, let $F\in (S/{\mathcal D})^*$ and let $f\in \{1,\dots,n+1\}^{S^{*}_F}$. It is clear that the map $g$ from ${\bf 2}^*$ to $(S/{\mathcal D})^*$ given by $g(a)=F$ is proper and continuous and so $f$ is the only component of a $g$-cohomomorphism from $S^*$ to ${\bf n+2}^*$. This and Theorem \ref{th1} imply that  the constructed assignment is surjective. \end{proof}

We now characterize the points of the space $\Lambda_n\lambda_n\Lambda_n(S)$. Let ${\mathcal T}_{n+1}$ denote the set of all transformations of the set $\{1,\dots, n+1\}$.

\begin{lemma}\label{lem201} There is a bijective correspondence between the points of the space $\Lambda_n\lambda_n\Lambda_n(S)$ and elements of the set
\begin{equation}\label{eq301}
\{(F,f,g)\colon F\in (S/{\mathcal D})^{*}, f\in\{1,\dots,n+1\}^{S^{*}_F}, g\in {\mathcal T}_{n+1}\}.
\end{equation}
\end{lemma}

\begin{proof} Let $h\in \Lambda_n\lambda_n\Lambda_n(S)$. That is, $h$ is a skew Boolean algebra morphism from $\lambda_n\Lambda_n(S)$ to ${\bf n+2}$. Then $\overline{h}$ is a Boolean algebra morphism from $(\lambda_n\Lambda_n(S))/{\mathcal D}$ to ${\bf 2}$. We have that $\overline{h}^{-1}(1)$ is an ultrafilter of $(\lambda_n\Lambda_n(S))/{\mathcal D}$ and so is equal to some $N_G$, where $G\in \Lambda_n(S)$, see \eqref{eq:nx}.  The only component of the cohomomorphism $\tilde{h}$ is a map from $(\lambda_n\Lambda_n(S))^{*}_{N_G}$ to $\{1,\dots,n+1\}$. This map defines $g\in{\mathcal T}_{n+1}$ where $N_{G,i}\mapsto g(i)$ for all $i\in\{1,\dots, n+1\}$. By Lemma \ref{lem20} we have that $G$ corresponds to a pair $(F,f)$. Hence $h$ corresponds to the triple $(F,f,g)$. That this correspondence is bijective by similar arguments as in the proof of Lemma \ref{lem20}.\end{proof}
 
In what follows, we often identify the points of the spaces $\Lambda_n(S)$ and $\Lambda_n\lambda_n\Lambda_n(S)$ with the elements of the sets \eqref{eq300} and \eqref{eq301}, respectively, via the constructions given in the proofs of Lemmas \ref{lem20} and \ref{lem201}.

To proceed, we need to establish how the action of maps of our interest looks like in this notation. We do this in the following three lemmas. By $id$ we denote the identity transformation of  the set $\{1,\dots, n+1\}$.

\begin{lemma}\label{lem202} 
$\epsilon_{\Lambda_n(S)}(F,f)=(F,f,id)$. \end{lemma}

\begin{proof} Let $\varphi\in \Lambda_n(S)$ and assume that $\epsilon_{\Lambda_n(S)}(\varphi)=\varphi'$. By Lemma \ref{lem:counit} 
$\varphi'(g)=g(\varphi)$, $g\in \lambda_n\Lambda_n(S)$.  Thus for $i=1,\dots, n+1$ we have
$$
g\in \varphi'^{-1}(i) \, \Leftrightarrow  \, \varphi \in g^{-1}(i)  \, \Leftrightarrow  \, g\in N_{\varphi, i}.
$$
Therefore $\varphi'^{-1}(i)=N_{\varphi, i}$. Thus also $\overline{\varphi'}^{-1}(1)=N_{\varphi}$. This and the constructions in the proofs of Lemmas \ref{lem20} and \ref{lem201} imply that if $\varphi$ corresponds to $(F,f)$ that $\varphi'$ corresponds to $(F,f,id)$, as required.
\end{proof}
 
\begin{lemma}\label{lem204} $\overline{\eta_S}^{-1}(N_{(F,f)})=F$  and the component $\widetilde{\eta_S}_{N_{(F,f)}}$ of $\widetilde{\eta_S}$ is given by
$$
x\mapsto N_{(F,f),f(x)}, x\in S^{*}_F.
$$
\end{lemma}

\begin{proof}
Let $\varphi\in \Lambda_n(S)$. Then for $a\in S$ and $i=1,\dots, n+1,$ applying the definition of $N_{\varphi,i}$ and \eqref{aux13}, we have
$$
a\in \eta_S^{-1}(N_{\varphi,i}) \,\Leftrightarrow \, \eta_S(a)\in N_{\varphi,i}\,\Leftrightarrow \, \eta_S(a)(\varphi)=i \, \Leftrightarrow \varphi(a)=i \,\Leftrightarrow \, a\in \varphi^{-1}(i).
$$
Thus $\eta_S^{-1}(N_{\varphi,i})=\varphi^{-1}(i)$.
If $\varphi$ corresponds to the pair $(F,f)$ via the construction in the proof of Lemma \ref{lem20} then the only component of $\widetilde{\varphi}$ is $f: S^*_F\to \{1,\dots, n+1\}$. We therefore have $\widetilde{\eta_S}^{-1}(N_{(F,f),i})=f^{-1}(i)$, and the statement follows.
\end{proof}

\begin{lemma}\label{lem203}
$\overline{\mu_S}^{-1} = \epsilon_{\Lambda_n(S)}^{**}$, that is, $\overline{\mu_S}^{-1}(N_x)=N_{\epsilon_{\Lambda_n(S)}(x)}$ and
the component $\widetilde{{\mu}_S}_{N_{(F,f)}}$ is given by 
$$
N_{(F,f,id),i}\mapsto N_{(F,f),i}.
$$
\end{lemma}
\begin{proof}
The statement follows from $\mu_S=\lambda_n(\epsilon_{\Lambda_n(S)})$ applying Lemma \ref{lem:aaa} and Lemma \ref{lem202}.
\end{proof}

We are now prepared to characterize the algebras for the monad $(T,\eta,\mu)$. By definition an {\em algebra of the monad} $(T,\eta,\mu)$ (or just a {\em $T$-algebra}) is a pair $(S,\gamma)$ with $S$  is a skew Boolean algebra and $\gamma:T(S)\to S$ is a morphism such that the following diagrams commute:

\begin{center}
\begin{tikzpicture}
\matrix (m) [matrix of math nodes, row sep=3.5em, column sep=2em, text height=1.5ex, text depth=0.25ex] 
{ T^2(S) & & T(S)  & & S & & T(S)\\
T(S) & & S & & & &  S\\};
\path[->] (m-1-1) edge node[above] {$T(\gamma)$} (m-1-3);
\path[->] (m-2-1) edge node[above] {$\gamma$} (m-2-3);

\path[->] (m-1-1) edge node[left] {$\mu_S$} (m-2-1);
\path[->] (m-1-3) edge node[right] {$\gamma$} (m-2-3);

\path[->](m-1-5) edge node[above]{$\eta_S$}(m-1-7);
\path[->](m-1-5) edge node[below]{$1_S$}(m-2-7);
\path[->](m-1-7) edge node[right]{$\gamma$}(m-2-7);
\end{tikzpicture}
\end{center}

A {\em morphism of $T$-algebras} $h:(S_1,\gamma)\to (S_2,\delta)$ is a morphism $h:S_1\to S_2$ such that the following diagram commutes:
\begin{center}
\begin{tikzpicture}
\matrix (m) [matrix of math nodes, row sep=3.5em, column sep=2em, text height=1.5ex, text depth=0.25ex] 
{ T(S_1) & & T(S_2) \\
S_1 & &S_2\\};
\path[->] (m-1-1) edge node[above] {$T(h)$} (m-1-3);
\path[->] (m-2-1) edge node[above] {$h$} (m-2-3);

\path[->] (m-1-1) edge node[left] {$\gamma$} (m-2-1);
\path[->] (m-1-3) edge node[right] {$\gamma$} (m-2-3);
\end{tikzpicture}
\end{center}

\begin{theorem}\label{th33}

\begin{enumerate}\label{th33}
\item A pair $(S,\gamma)$ is an algebra for the monad $(T,\eta,\mu)$ if and only if $S=\lambda_n(X)$ for some Boolean space $X$ and $\gamma=\lambda_n\epsilon_X$.
\item A map $h:\lambda_n(X_1)\to \lambda_n(X_2)$ is a morphism of $T$-algebras if and only if $h=\lambda_n(f)$ for some morphism $f:X_1\to X_2$.
\end{enumerate}
\end{theorem}

\begin{proof} 
The equality $1_S=\gamma\eta_S$ implies the equality $1_{(S/D)^{*}}={\overline{\eta_S}}^{-1}{\overline{\gamma}}^{-1}$. Applying Theorem \ref{th1} we have the following commuting diagrams:
\begin{center}
\begin{tikzpicture}
\matrix (m) [matrix of math nodes, row sep=3.5em, column sep=2em, text height=1.5ex, text depth=0.25ex] 
{ S^{*} & & (\lambda_n\Lambda_n(S))^{*}  & & (S/{\mathcal D})^{*} & & (\lambda_n\Lambda_n(S)/{\mathcal D})^*\\
 & & S^{*} & & & &  (S/{\mathcal D})^{*}\\};
\path[->] (m-1-1) edge node[above] {$\widetilde{\eta_S}$} (m-1-3);
\path[->] (m-1-3) edge node[right]{$\widetilde{\gamma}$} (m-2-3);

\path[->] (m-1-1) edge node[left] {$1_{S^{*}}$} (m-2-3);

\path[<-](m-1-5) edge node[above]{$\overline{\eta_S}^{-1}$}(m-1-7);
\path[<-](m-1-5) edge node[below]{$1_{(S/{\mathcal D)^{*}}}\,\,\,\,\,\,$}(m-2-7);
\path[<-](m-1-7) edge node[right]{$\overline{\gamma}^{-1}$}(m-2-7);
\end{tikzpicture}
\end{center}

 Let $F\in(S/\D)^{*}$. Assume $\overline{\gamma}^{-1}(F)=N_{(G,f)}$. By Lemma \ref{lem204} we have $\overline{\eta_S}^{-1}(N_{(G,f)})=G$. Since the second diagram above commutes we obtain that $F=G$. We fix $f=f_F: S^{*}_F\to \{1,\dots, n+1\}$ such that $\overline{\gamma}^{-1}(F)=N_{(F,f)}$. Now, since the first diagram above commutes and applying Lemma \ref{lem204}, we see that for $x\in S^{*}_F$ it holds
\begin{equation}\label{eq15}
x\stackrel{\widetilde{\eta_S}_{N_{(F,f)}}}{\xrightarrow{\hspace{1cm}}} N_{(F,f), f(x)} \stackrel{\widetilde{\gamma_F}}{\xrightarrow{\hspace{1cm}}} x.
\end{equation}
The latter shows that the map $\widetilde{\eta_S}_{N_{(F,f)}}$ and  also the restriction of the map map $\widetilde{\gamma_F}$ to the image of $\widetilde{\eta_S}_{N_{(F,f)}}$ are  injective. Injectivity of $\widetilde{\eta_S}_{N_{(F,f)}}$ implies that $f$ must be injective, too.

By definition and applying Theorem \ref{th1} we have the following diagrams (where $\Lambda_n(\gamma)^{**}(N_x)=N_{\Lambda_n(\gamma)(x)}$):

\begin{center}
\begin{tikzpicture}
\matrix (m) [matrix of math nodes, row sep=3.5em, column sep=0.6em, text height=2ex, text depth=0.25ex] 
{ (\lambda_n\Lambda_n\lambda_n\Lambda_n(S))^{*} & & & &(\lambda_n\Lambda_n(S))^{*} & (\lambda_n\Lambda_n\lambda_n\Lambda_n(S)/{\mathcal D})^*& & & & (\lambda_n\Lambda_n(S)/{\mathcal D})^*\\
(\lambda_n\Lambda_n(S))^{*} & & & & S^{*}&  (\lambda_n\Lambda_n(S)/{\mathcal D})^*& & & & (S/{\mathcal D})^{*}\\};
\path[->] (m-1-1) edge node[above] {$\widetilde{\lambda_n\Lambda_n(\gamma)}$} (m-1-5);
\path[->] (m-2-1) edge node[above] {$\widetilde{\gamma}$} (m-2-5);
\path[->] (m-1-1) edge node[left] {$\widetilde{\mu_S}$} (m-2-1);
\path[->] (m-1-5) edge node[right] {$\widetilde{\gamma}$} (m-2-5);

\path[<-] (m-1-6) edge node[above] {$\Lambda_n(\gamma)^{**}$} (m-1-10);
\path[<-] (m-2-6) edge node[above] {$\overline{\gamma}^{-1}$} (m-2-10);
\path[<-] (m-1-6) edge node[left] {$\epsilon_{\Lambda_n(S)}^{**}$} (m-2-6);
\path[<-] (m-1-10) edge node[right] {$\overline{\gamma}^{-1}$} (m-2-10);

\end{tikzpicture}
\end{center}

Our goal now is to describe the action of $\Lambda_n(\gamma)$ and $\lambda_n\Lambda_n(\gamma)$ if the elements of the spaces are encoded as is described in the proofs  of Lemmas \ref{lem20} and  \ref{lem201}. Let $g\in\Lambda_n(S)$. Its corresponding pair is $(F,\hat{g})$, where $F=\overline{g}^{-1}(1)$ and $\hat{g}$ is the only component of $\widetilde{g}$. Applying $\Lambda_n(\gamma)(g)=g\gamma$ and the construction in the proof of Lemma \ref{lem201} it follows that  we have
\begin{equation}\label{eq310}
(F,\hat{g})\stackrel{\Lambda_n(\gamma)}{\xrightarrow{\hspace{1cm}}} (F,f,\hat{g}\widetilde{\gamma}_F).
\end{equation}

This and Lemma \ref{lem:aaa} yield
\begin{equation}\label{eq311}
N_{(F,f,\hat{g}\widetilde{\gamma_F}),i}\stackrel{\widetilde{\lambda_n\Lambda_n(\gamma)}_{(F,\hat{g})}}{\xrightarrow{\hspace{1.5cm}}} N_{(F,\hat{g}),i}
\end{equation}
for each $i\in\{1,\dots, n+1\}$. It follows that for $F\in (S/{\mathcal D})^{*}$ we have

$$
F \stackrel{\overline{\gamma}^{-1}} {\xrightarrow{\hspace{1cm}}}N_{(F,f)} \stackrel{\Lambda_n(\gamma)^{**}}{\xrightarrow{\hspace{1cm}}} N_{(F,f,f\widetilde{\gamma}_F)}.
$$

On the other hand,  observe that $(\lambda_n\Lambda_n(S)/{\mathcal D})^*$  is isomorphic to $\Lambda_n(S)$ and $(\lambda_n\Lambda_n\lambda_n\Lambda_n(S)/{\mathcal D})^*$  is isomorphic to $\Lambda_n\lambda_n\Lambda_n(S)$ via the map $N_x\mapsto x$. Applying  Lemma \ref{lem202}, for $F\in (S/{\mathcal D})^{*}$ we have
$$
F \stackrel{\overline{\gamma}^{-1}}{\xrightarrow{\hspace{1cm}}} N_{(F,f)} \stackrel{\epsilon_{\Lambda_n(S)}^{**}}{\xrightarrow{\hspace{1cm}}} N_{(F,f,id)}.
$$

Since the two expressions above must be equal, it follows that $f$ is a bijection and
$\widetilde{\gamma}_F=f^{-1}$.
Hence $\vert S^{*}_F\vert = n+1$. Therefore, all the stalks of $S^*$ are $(n+1)$-element.

Let $F\in (S/\D)^{*}$ be fixed. Since $\vert S^{*}_F\vert = n+1$ and $f$ is a bijection, we can enumerate the germs of $S^{*}_F$ so that $S^{*}_F=\{F_{(1)},\dots,F_{(n+1)}\}$ and $f$ maps each $F_{(i)}$ to $i$. In this notation we have
\begin{equation}\label{aux16}\widetilde{\gamma}_F(N_{(F,f),i})=F_{(i)}.
\end{equation}
Let $i\in\{1,\dots, n+1\}$. We define the {\em $i$-th layer} of the space $S^{*}$ as the set $S^{*}_{(i)}=\{F_{(i)}\colon F\in (S/\D)^{*}\}$.

We proceed to show that $S\simeq \lambda_n((S/{\mathcal D})^{*})$. It is convenient to work with the isomorphic copy $S^{**}$ of $S$.  Let $s\in S$. Observe that $\widetilde{\eta_S}(M(s)\cap S^{*}_{(i)})$ is the $i$-th layer of $M(\eta_S(s))$. Since the projection of the latter set is compact-open then its inverse image under $\overline{\gamma}$ is compact-open, too. 
It follows that the map $F_{(i)}\mapsto N_{F,i}$,  where $F\in (S/{\mathcal D})^{*}$ and $i\in\{1,\dots, n+1\}$, induces an isomorphism between $S^{**}$ and $\lambda_n((S/{\mathcal D})^{*})^{**}$.

Let $X=(S/{\mathcal D})^{*}$. To establish the equality $\gamma=\lambda_n\epsilon_X$ we apply arguments similar to those in the proofs of Lemma \ref{lem202} and Lemma \ref{lem203} to observe that the action of  $\widetilde{\lambda_n\epsilon_X}_F$ coincides  with the action of $\widetilde{\gamma}_F$ given in \eqref{aux16}.

We are left to prove the claim about morphisms. Assume first that $f:X_2\to X_1$ is a morphism of Boolean spaces. It is straightforward to verify that  $\lambda_n(f):\lambda_n(X_1)\to \lambda_n(X_2)$ is a morphism of $T$-algebras from $(T(\lambda_n(X_1)),\lambda_n\epsilon_{X_1})$ to $(T(\lambda_n(X_2)),\lambda_n\epsilon_{X_2})$.  

We now assume that $h:\lambda_n(X_1)\to \lambda_n(X_2)$ is a morphism of $T$-algebras from $(T(\lambda_n(X_1)),\lambda_n\epsilon_{X_1})$ to $(T(\lambda_n(X_2)),\lambda_n\epsilon_{X_2})$. Since $(\lambda_n(X_i)/{\mathcal D})^{*}\simeq X_i$, $i=1,2$, then $\overline{h}^{-1}$ induces a morphism, $\hat{h}$, from  $X_2$ to $X_1$. We show that $h=\lambda_n(\hat{h})$. Let $(F,f)\in\Lambda_n(S_2)$, where $F\in X_2$ and $f\in \{1,\dots,n+1\}^{(\lambda_n(X_2))^{*}_F}$. Then 
$$
(F,f)\stackrel{\Lambda_n(h)}{\xrightarrow{\hspace{1cm}}} (\overline{h}^{-1}(F),f\widetilde{h}_F)
$$
and therefore
\begin{equation}\label{bb}
N_{(\overline{h}^{-1}(F),f\widetilde{h}_F),i} \stackrel{\widetilde{\lambda_n\Lambda_n(h)}_{N_{(F,f)}}}{\xrightarrow{\hspace{1.5cm}}} N_{(F,f),i}
\end{equation}
for all $i=1,\dots, n+1$.

Let $X$ be a Boolean space. Similarly as it was done in Lemma \ref{lem201} we establish a bijection between the points of $\Lambda_n\lambda_n(X)$ and pairs $(F,f)$, $F\in X$, $f\in {\mathcal T}_{n+1}$. If $g\in \Lambda_n\lambda_n(X)$, then $\overline{g}: X^{*}\to {\bf 2}$ determines $F$ and the only component of $\widetilde{g}$ determines $f$. By Lemmas \ref{lem:aaa} 
in this notation we have 
\begin{equation}\label{cc} N_{(F,id),i}\stackrel{\widetilde{\lambda_n\epsilon_{X_2}}_F}{\xrightarrow{\hspace{1cm}}}N_{(F,i)}.
\end{equation}
From \eqref{bb},\eqref{cc} and the commutative diagram
\begin{center}
\begin{tikzpicture}
\matrix (m) [matrix of math nodes, row sep=3.5em, column sep=2em, text height=1.5ex, text depth=0.25ex] 
{ \lambda_n\Lambda_n\lambda_n(X_1) & & \lambda_n\Lambda_n\lambda_n(X_2) \\
\lambda_n(X_1) & &\lambda_n(X_2)\\};
\path[->] (m-1-1) edge node[above] {$\lambda_n\Lambda_n(h)$} (m-1-3);
\path[->] (m-2-1) edge node[above] {$h$} (m-2-3);

\path[->] (m-1-1) edge node[left] {$\lambda_n\epsilon_{X_1}$} (m-2-1);
\path[->] (m-1-3) edge node[right] {$\lambda_n\epsilon_{X_2}$} (m-2-3);
\end{tikzpicture}
\end{center}
we see that $\widetilde{\lambda_n\epsilon_{X_1}}$ is defined on all  $ N_{(\overline{h}^{-1}(F), \widetilde{h}_F),i}$,  $F\in X_2$, $1\leq i\leq n+1$, and it must be $\widetilde{h}_F=id$.  It follows that the action of $\widetilde{h}_F$ is given by $N_{\overline{h}^{-1}(F),i}\mapsto N_{F,i}$ and thus $h=\lambda_n(\hat{h})$. The proof is complete.
\end{proof}

Let $\lambda_n(\LCBS)$ be the category whose objects are $\lambda_n(X)$, where $X$ is a Boolean space, and whose arrows are $\lambda_n(f)$, $f$ is a morphism of Boolean spaces.

\begin{corollary} The Eilenberg-Moore category of the monad $(T,\eta,\mu)$ is isomorphic to the category $\lambda_n(\LCBS)$. Consequently, the adjunction $\Lambda_n\dashv \lambda_n$ is monadic for every $n\geq 0$.
\end{corollary}

\begin{proof}
The first statement follows from Theorem \ref{th33}. The second statement holds because the category $\lambda_n(\LCBS)$ is obviously isomorphic to the category $\LCBS^{op}$.
\end{proof}

\begin{corollary}
The category $\lambda_n(\LCBS)$ is a reflective subcategory of the the category $\LSBA$.
The reflector is given by the functor $\Lambda_n\lambda_n$.
\end{corollary}

\end{document}